\documentclass[a4paper,12pt]{article}

\usepackage[utf8x]{inputenc}
\usepackage{enumerate}
\usepackage{amsmath}
\usepackage{textcomp}

\title{Generalization of Strongly Clean Rings}
\author{Abhay K. Singh\\Department of Applied Mathematics\\Indian School of Mines, Dhanbad\\India}
\date{}
\begin{document}
\maketitle
\begin{center}\textbf{Abstract}\end{center}
In this paper, strongly clean ring defined by W. K. Nicholson [6] in 1999 has been generalized to n-strongly clean, ∑-strongly clean and with the help of example it has been
shown that there exists a ring, which is n-strongly clean and ∑-strongly clean but not strongly clean. It has been shown that for a commutative ring R formal power series
R[( x )] of R is n-strongly clean if and only if R is n- strongly clean. We also discussed the structure of homomorphic image of n- strongly clean and direct product of n- strongly
clean rings. It has also been shown that for any commutative ring R, the polynomial ring R ( x ) is not ∑-strongly clean ring.
\newline\newline
\textbf{Keywords-} Strongly clean rings, Clean rings, $\sum$-Clean, n-Clean rings, Group rings, Exchange rings.
\newline\newline
(2010) Mathematics subject Classification: 16D70, 16P70, 16E50
\newline\newline
\textbf{1:Introduction}
\newline\newline 
\indent Throughout this paper, we denote R an associative ring with unity. Let m
Z$_{p}$ = {$\frac{m}{n}$, m.n ∈ Z and gcd(p,n) = 1} , where p is prime number and Z is ring of integer is a subring of the ring of rational numbers Q. The Jacobson radical of R is denoted by J(R)
and Id(R), U(R) denote the idempotent element in R, unit element in R respectively.
\newline\newline
\textit{Definition 1.1:} An element x of a ring R is strongly clean if x = e + u and eu = ue ,where e$^{2}$ = e and u is unit in R. 
\newline\newline
\textit{Definition 1.2:} Let n be any positive integer. An element x of R is called n-strongly clean ring if x= e+u$_{1}$+u$_{2}$+ . . .  +u$_{n}$ such that $eu_{1} = u_{1}e$, $eu_{2} = u_{2}e$,
 . . . , $eu_{n} = u_{n}e$ where e is idempotent and u$_{1}$, u$_{2}$, . . . , u$_{n}$ are units in R. A ring R is called n-strongly clean if every element of R is n-strongly clean.
\newline\newline
\textit{Definition 1.3:}  Let G be a cyclic group of order 3 then group ring Z$_{p}$G is 2-strongly clean but not strongly clean for any prime p$≠$2.
\newline\newline
\textit{Proof:}  Let G = {1, a, a$^{2}$}, with a$^{3}$ =1 then,\newline
\indent Z$_{p}$G = $\{$ $\frac{m_{1}}{n_{1}}$ + $\frac{m_{2}}{n_{2}}$ + $\frac{m_{3}}{n_{3}}a^{2}/m_{i},n_{i}$ $\in$ Z and p does not divide $n_{i}\}$,
\newline where i = 1,2,3 is a commutative. \newline
\indent Y. Yuanqing, [Theorem 3.1, [5]] showed that Hence Z$_{p}$G is semiclean. Hence Z$_{p}$G= p + u , where p is periodic and u is unit in Z$_{p}$G.
\newline\indent By [Lemma 5.1, [5]] periodic elements are clean, p can be written as\begin{center}p = e+v,\end{center}
\indent Where e is an idempotent element in Z$_{p}$G and v is unit in ZpG 
\newline\indent So we obtain \begin{center}Z$_{p}$G = e + v + u ,\end{center}
Hence Z$_{p}$G is 2- strongly clean.\newline\indent
Han and Nicholson [3] showed that Z$_{7}$G is not clean. Hence, 2-strongly clean ring need not be strongly clean.
\newline\newline
\textit{Proposition 1.4:}Let n be a positive integer. Then following conditions are equivalent:
\begin{enumerate}[(i)]
 \item Homomorphic image of n-strongly clean ring is n-strongly clean.
 \item A direct product R = $\Pi_{\alpha}$R$_{\alpha}$of ring (R$_{\alpha}$) is n-strongly clean if and only if each R$_{\alpha}$ is n–strongly clean
\end{enumerate}
$ $\newline
\textit{Proof:}  (i) is obvious, so we need to prove only (ii)
\newline\newline
\indent Suppose each R$_{\alpha}$ is an n–strongly clean ring and for each $\alpha$ , write
x$_{\alpha}$ = e$_{\alpha}$ + u$_{\alpha}^{1}$ + . . . + u$_{\alpha}^{n}$ , where u$_{\alpha}^{i}$ $\in$ U(R$_{\alpha}$)(1$\leq$ i $\leq$ n ) and e$_{\alpha}$ $\in$ Id(R ) such that every 
e$_{\alpha}$ commutes with every u$_{\alpha}$. Then x = e + u$_{1}$ + u$_{2}$ + . . . + u$_{n}$, where u$_{\alpha}^{i}$ $\in$ U(R$_{\alpha}$)(1$\leq$ i $\leq$ n ) and e$_{\alpha}$ $\in$ 
Id(R ) is an n- strongly clean ring.
\newline\newline
\textit{Proposition 1.5}: Let R be a commutative ring and n be a positive integer, then R [(x)] is n-strongly clean if and only if R is n- strongly clean.
\newline\newline
\textit{Proof:}  If R[x] is n-strongly clean ring then from proposition 4.3 (i), R=R[x]|(x) is n-strongly clean ring.\newline\newline
Conversely, suppose that R is n-strongly clean, then let
\begin{center} $f = \sum \limits_{i=0}^{\infty}r_{i}x^{i} \in R[x]$ write,\end{center}
r$_{o}$ = e + u$_{1}$ + u$_{2}$ + . . . + u$_{n}$, eu$_{1}$ = u$_{1}$e, eu$_{2}$ = u$_{2}$e, . . ., eu$_{n}$ = u$_{n}$e where e $\in$ Id(R) u$_{1}$, u$_{2}$, u$_{3}$,
. . ., u$_{n}$ $\in$U(R) then,\newline
f = e + (u$_{1}$ + r$_{1}$x + r$_{2}$x$^{2}$+ . . .) + u$_{2}$ + . . .+ u$_{n}$ where e $\in$ Id(R) $\subseteq$ Id(R[x]). \newline
So, u$_{1}$ + r$_{1}$x + r$_{2}$x$^{2}$+ . . . $\in$ U[R(x)], \newline
u$_{i}\in$U(R)$\subseteq$U[R(x)] (1$\leq$i$\leq$n).
\newline\newline
Hence R[x] is n-strongly clean ring.
\newline\newline
\textbf{Proposition 1.6:} Let n be a positive integer and I is an ideal of R, which is contained in J(R). If R/I is n- strongly clean and idempotent can be lifted modulo I, then R is n-
strongly clean.
\newline\newline
\textbf{Proof:} Let a$\in$R and \={a} = \={e} + \={u$_{1}$} + \={u$_{2}$} + . . . + \={u$_{n}$} where e$^{2}$ - e $\in$ I and \={u$_{i}$v$_{i}$} = \={v$_{i}$u$_{i}$} = \={1}
with v$_{i}$ R(i = 1,2,3, . . . ,n). Since idempotent can be lifted modulo I, we may assume that e$^{2}$ = e$\in$R. Also we can assume that u$_{i}$(1$\leq$i$\leq$n)
is unit in R because u$_{i}$v$_{i}$ = 1 + r $_{i}$, v$_{i}$u$_{i}$ = 1 + t$_{i}$ for some r$_{i}$, t$_{i}$ $\in$I$\subseteq$J(R), but J(R)=\{r$\in$R/r+a is a unit for every
unit in a$\in$R\}. So R is a n-strongly clean.
\newline\newline\indent
Camillo and Yu [2] in 1994 showed that R is clean if and only if every element of R is sum of unit and square root of 1.
\newline\newline
\textbf{Corollary 1.7:} Let n be a positive integer and R be a commutative ring. Then R is strongly clean if and only if every element of R is sum of unit and a square root of 1.
\newline\newline
We extended this result to n-strongly clean rings.
\newline\newline
\textbf{Proposition 1.8:} Let R be a commutative ring in which 2 is invertible and n be a positive integer .Then R is n- strongly clean if and only if R is sum of n-units and square root of 1.
\newline\newline
\textbf{Proof:} Suppose R is n- strongly clean and x ∈ R , then \newline 
$\frac{x+1}{2}$ = e + u$_{1}$ + u$_{2}$ + . . .  + u$_{n}$,\newline
where e$^{2}$ = e, eu$_{1}$ = u$_{1}$e, eu$_{2}$ = u$_{2}$e, . . ., eu$_{n}$ = u$_{n}$e and\newline
thus x = (2e-1) + 2u$_{1}$ + 2u$_{2}$ + . . . + 2u$_{n}$ with (2e-1)$^{2}$ = 1 and 2u$_{1}$,2u$_{2}$,. . . ,2u$_{n}$ = U(R).\newline
Conversely, if x$\in$R then 2x-1 = f + u$_{1}$ + u$_{2}$ + . . .  + u$_{n}$ where f$^{2}$=f and u$_{1}$,u$_{2}$,. . . ,u$_{n}$ = U(R).\newline
Thus x = $\frac{1}{2}(f+1)$ + $\frac{1}{2}$(u$_{1}$) + $\frac{1}{2}$(u$_{2}$) + . . . + $\frac{1}{2}$(u$_{n}$) $\in$ U(R).
Now, $[\frac{1}{2}(f+1)]^{2}$ = $[\frac{1}{2}(f+1)]$ $\in$ Id(R) and $\frac{1}{2}$u$_{1}$,. . .,$\frac{1}{2}$u$_{n}$ $\in$U(R) and due to commutibility of R, x is\newline
n-strongly clean element.
\newline\indent
For any positive integer n, we let U$_{n}$(R ) denote the set of the element of R that can be written as no sum of more than n units of R. A ring R is called generated by its
units if every elements of R is sum of finitely many units of R. Camillo and Yu[2] in 1994 showed that if R is clean ring then R = U$_{2}$(R ) . Now if R is strongly clean
then R = U$_{2}$(R ) , and with above proposition we have following result for n- strongly clean ring.
\newline\newline
\textbf{Corollary 1.19:} Let R be a ring in which 2 is invertible and n be a positive integer .If R is n-strongly clean then R = U$_{n +1}$(R ) .
\newline\newline
\textbf{Theorem 1.1.0:} Let n be any positive integer and e is idempotent of R. If e Re and (1 − e)R(1 − e) are both n-strongly clean then R is also n- strongly clean.
\newline
\textbf{Proof:} We use \={e} to denote the 1-e and apply Pierce decomposition for the ring R, i.e.,\newline
R = ere + er\={e} + \={e}Re + \={e}r\={e}. Let x $\in$R. The we can write\begin{center}x = a+b+c+d,\end{center}
Where a,b,c,d belongs to ere,\={e}re,\={e}r\={e}, respectively. By our assumption a is n-strongly clean. There exist f,u$_{1}$,u$_{2}$,. . .,u$_{n}$, fu$_{1}$ = 
u$_{1}$f, fu$_{2}$ = u$_{2}$f . . ., fu$_{n}$ = u$_{n}$f where f$^{2}$ = f $\in$ eRe and u$_{i}$(1$\leq$i$\leq$n) are units in eRe with inverse u$_{i}^{-1}$. Thus,
d-cu$^{-1}$b$\in$\={e}R\={e}. Therefore, \begin{center}g-cu$^{-1}$b = g + v$_{1}$ + v$_{2}$ + . . . + V$_{n}$ \end{center}
Where g$^{2}$=g $\in$\={e}R\={e}, gv$_{1}$ = v$_{1}$g, gv$_{2}$ = v$_{2}$g, . . ., g$v_{n}$ = v$_{n}$g and v$_{i}$(1$\leq$i$\leq$n) are units in \={e}R\={e} with
inverse v$_{i}^{-1}$.\newline\newline
Hence,\newline\newline
\indent x = a + b + c + d\newline
\indent   = (f + u$_{1}$ + . . . + u$_{n}$) + b + c + (g + v$_{1}$ + . . . + v$_{n}$).\newline
\indent   = (f + g) + (u$_{1}$ + b + c + v$_{1}$ + c$u_{1}^{-1}$b) + (u$_{2}$ + v$_{2}$) + . . . + (u$_{n}$ + v$_{n}$).\newline
We will prove that f + g is an idempotent element of R and u$_{1}$ + b + c + v$_{1}$ + c$u_{1}^{-1}$b + v$_{i}$(2$\leq$i$\leq$n) are units in R.\newline
Now, the following fact is clear.\newline
(u$_{1}$ + b + c + v$_{1}$ + c$u_{1}^{-1}$b)(u$_{1}^{-1}$ + u$_{1}$bv$_{1}^{-1}$cu$_{1}^{-1}$ - u$_{1}^{-1}$bv$_{1}^{-1}$ - v$_{1}^{-1}$cu$_{1}^{-1}$ + v$_{1}^{-1}$) = 1.\newline
(u$_{i}$ + v$_{i}$)(u$_{i}^{-1}$ + v$_{i}^{-1}$) = e +(1-2) = 1 (2$\leq$i$\leq$n)\newline
We have, x = e$^{1}$ + p$_{1}$ + p$_{2}$ + . . . + p$_{n}$, where e$^{1}$ = f + , p$_{1}$ = u$_{1}$ + b + c + v$_{1}$ + cu$_{1}^{-1}b$, and p$_{i}$ = u$_{i}$ +v$_{i}$
(2$\leq$i$\leq$n). Hence e$^{1}$ is an idempotent element of R, p$_{1}$,.........p$_{n}$ are units in R and idempotent commute with every units. This proves that x is n-strongly clean.
\newline\newline
\textbf{Corollary 1.11:} Let n be positive integer. If 1 = e$_{1}$ + e$_{2}$ + e$_{3}$ + . . . + e$_{n}$ in a ring R, where e$_{i}$ orthogonal idempotent elements and each e$_{i}$Re$_{i}$ 
is n-strongly clean, then R is n-strongly clean.
\newline
\textbf{Proof:} Followed by Theorem 1.10
\newline\newline
\textbf{Corollary 1.13:} Let n be a positive integer. If M = M$_{1}\oplus$M$_{2}\oplus$ . . . $\oplus$M$_{m}$ are modules and End(M$_{i}$) is n-strongly clean for each i, then End(M) is 
n- strongly clean.
\newline
\textbf{Proof}:  Followed by Theorem 1.10 and Corollary 1.11
\newline\newline
\textbf{Definition 1.13} An element x of a ring R is called ∑-strongly clean if there exist a positive integer n such that x is n- strongly clean. A ring R is called a ∑-strongly clean if
every element of R is a sum of idempotent and a finite number of units in R such that idempotent commutes with every unit.\newline\newline\indent
Clearly strongly clean rings are ∑-strongly clean rings. Let R = Z , the ring of integer. Then R is a ∑-strongly clean ring, but it is not n- strongly clean for any positive
integer n. Hence ∑-strongly clean rings need not be n- strongly clean for any positive integer n. Hence the classes of strongly clean rings are a proper subset of ∑-strongly
clean rings.\newline\newline\indent
It has been proved by Xiao and Tong [1] in 2005, if G is a cyclic group of order 3 then group ring Z$_{2}$G is clean. We can extend this result to strongly clean ring.
\newline\newline
\textbf{Corollary 1.14:} If G is cyclic group of order 3 then group ring Z2G is strongly clean.\newline\newline
By corollary 1.14 and Example 1.3, the following result is immediate.
\newline\newline
\textbf{Corollary 1.15:} f G is a cyclic group of order 3, then group ring ZpG is ∑-strongly
clean ring for any prime p$≠$2\newline\newline
It is obvious that the polynomial ring over a field is not strongly clean. Anderson and Camillo [7] in 2002 showed that every polynomial ring over a commutative ring is
not clean hence not strongly clean. The following example shows that every polynomial over commutative ring is not ∑-strongly clean.
\newline\newline
\textbf{Example 1.16:} If R is any commutative ring, then the polynomial ring R[ x ] is not ∑-strongly clean.
\newline
\textbf{Proof:} Note that Id(R[ x ]) = Id(R ) and,\newline\newline
U(R[x] = \{r$_{o}$ + r$_{1}$x + . . . + r$_{n}$x$^{n}$/r$_{o}$ $\in$ U(R), r$_{i}$ $\in$ $\sqrt{0}$, i = 1,2,...,n\} where $\sqrt{0}$ denote the nil radical of R.
\newline\newline
If x is ∑-strongly clean, we may let\newline\newline
\indent x = f + (u$_{1}$ + r$_{1}$x + . . .) + (u$_{2}$ + r$_{2}$x + . . .) + ... + (u$_{n}$ + r$_{n}$x + . . .),\newline
where f$^{2}$ = f $\in$ R, u$_{1}$, u$_{2}$, ..., u$_{n}$ $\in$ U(R) and r$_{1}$,r$_{2}$,...,r$_{n}$ $\in$ $\sqrt{0}\subseteq$J(R),(1$\leq$i$\leq$). Then
$\sum\limits_{i=1}^{n}$ = 1 $\in$ J(R), which is a contradiction. Thus R[ x ] is not is ∑-strongly clean.
\newline\newline
\textbf{Proposition 1.17:} The following conditions are equivalent:
\begin{enumerate}[(1)]
 \item A homomorphic image of a ∑- strongly clean ring is ∑-strongly clean ring.
 \item A direct product R = $\Pi$R$_{\beta}$of ring (R$_{\beta}$ ) is ∑-strongly clean ring, then each R$_{\beta}$ is ∑-strongly clean ring.
\end{enumerate}
\textbf{Proof:} Clear from the proof of proposition 1.3
\newline\newline
\textbf{Proposition 1.18:} Let n be a positive integer and I be an ideal of R that is contained in J(R). If R/I is ∑-strongly clean and idempotent can be lifted 
modulo I, then R is ∑-strongly clean.
\newline
\textbf{Proof:} Followed from the proof of proposition 1.7.
\newline\newline
\textbf{Proposition 1.19:} Let n be positive integer and R be a commutative ring in which 2 is invertible .Then R is ∑-strongly clean if and only if every element 
of R is sum of n units and a square root of 1.
\newline\newline
\textbf{Proof:} Followed by the proof of Proposition 1.7.
\newline\newline
\textbf{Corollary 1.21:} Let R be a ring in which 2 is invertible and n be positive integer. If R is
∑-strongly clean, then R = U$_{n+1}$(R ).\newline
\textbf{Proof:} Followed from Corollary 1.9.
\newline\newline
\textbf{References}\newline\newline
1. Xiao Guangshi, Tong Wenting (2005), n-clean rings and weakly unit stable range rings , Comm. Algebra, 33, 1501-1517.\newline\newline
2. Camillo,V.P.,Yu,H.P. (1994) Exchange ring, units and idempotents, Comm. Algebra 22, 4737-4749.\newline\newline
3. Han,J, Nicolson,W.K. (2001) Extension of clean rings, Comm. Algebra 29,2589-2595.\newline\newline
4. Nicolson,W.K.(1977) Lifting’s idempotent and Exchange rings, Trans. Amer.Math.Soc., 229, 269-278.\newline\newline
5. Yuanqing Ye (2003), Semiclean rings, Comm. Algebra 31 (11) 5609-56256. Nicolson, W.K. (1999), Strongly clean rings and fitting lemma, Comm. Algebra, 27,3583-3592.\newline\newline
7. Anderson, D.D., Camillo, V.P., (2002) Commutative rings whose elements are some of unit and idempotent, Comm. Algebra, 30(7), 3327-3336\newline\newline
8:  Abhay Kumar Singh and B.M. Pandeya, A note on Generalization of Semiclean Rings International Journal of Algebra, Vol. 5, 2011, no. 21, 1039 – 1047.\newline\newline
\end{document}